\documentclass[a4paper,12pt]{article}
\usepackage{times, url}
\usepackage{graphicx,amsmath,amssymb,amsthm,amsfonts,bm,float,cite,soul,color}
\newtheorem{theorem}{Theorem}[section]

\newtheorem{lemma}{Lemma}[section]

\usepackage[none]{hyphenat}

\usepackage[center]{caption}
\usepackage[bottom,symbol]{footmisc}
\DeclareFontEncoding{TS1}{}{}
\DeclareFontSubstitution{TS1}{cmr}{m}{n}
\DeclareTextSymbol{\tcrp}{TS1}{'251}
\DeclareTextSymbolDefault{\tcrp}{TS1}

\usepackage[hidelinks]{hyperref} 
\usepackage{orcidlink}

\hypersetup{pdftitle={Your paper title here in sentence case}, pdfauthor={First Author, Second Author}}

\begin{document}
\setcounter{page}{1}

\thispagestyle{empty} 

\begin{center}
{\LARGE \bf  Tribonacci-Lucas numbers that are concatenations of two Tribonacci-Lucas numbers} 
\vspace{6mm}

{\Large \bf Kaggwa Prosper$^1$}
\vspace{4mm}

$^1$ Department of Mathematics, Makerere University \\
P.O. Box 7062, Kampala, Uganda \\
e-mail: \url{kaggwaprosper58@gmail.com}
\vspace{2mm}

\end{center}
{\bf Abstract:} The sequence of Tribonacci-Lucas numbers  $\{S_n\}_{n\ge 0}$ is defined by the recurrence relation $S_{n} = S_{n-1} + S_{n-2} + S_{n-3}$, with initial conditions $S_0 =S_2= 3$ and $S_1 = 1$. In this paper, we determine all Tribonacci-Lucas numbers that can be expressed as concatenations of two terms of the same sequence. We use Baker's theory for nonzero linear bounds for linear forms in logarithms of algebraic numbers and reduction methods based on the theory of continued fractions. 
\vspace{5mm}

\section{Introduction}
For $S_0 =S_2=3$ and $S_1=1,$ $S_{n+3} = S_{n+2} + S_{n+1} + S_n$ denotes the $n^{th}$ term of the Tribonacci-Lucas sequence for $n \ge 0$. The first few terms of this sequence are given by
\begin{equation*}
	(S_n)_{n \geq 0} = 3, 1, 3, 7, 11, 21, 39, 71, 131, 241, 443, 815, 1499, 2757, 5071, 9327, 17155, 31553, \ldots.
\end{equation*}

This sequence shares the same recurrence relation with the Tribonacci sequence only with different initial terms. Growing research has examined which terms of a sequence can be written as concatenations of other terms within the same sequence notably in \cite{alan, erduvan}, where the authors studied these existences in the $k$-Fibonacci sequence and Lucas sequence respectively. However, no such characterization is known for the Tribonacci-Lucas sequence, which motivates the present study. In this paper, we extend this work to the Tribonacci-Lucas sequence by determining solutions to the following Diophantine equation. 
\begin{equation}\label{main}
	S_k=S_m \cdot 10^d +S_n.
\end{equation} in non negative integers $k,m,n,d$ where $d$ is the number of digits of $S_n$. We state the main result of this paper in the following Theorem.
\begin{theorem}\label{thm1x}
	The only Tribonacci-Lucas numbers that are concatenations of two Tribonacci-Lucas numbers are $ S_4=11=\overline{S_1S_1} , S_7=71=\overline{S_3S_1}$.
\end{theorem}

This result gives a continuation of the results presented  in \cite{erduvan, alan}. The method of proof of our main results makes use of Baker's theory for nonzero lower bounds for linear forms in logarithms of algebraic numbers, and reduction methods involving the theory of continued fractions. All remaining cases were verified computationally using SageMath.

\section{Preliminary results}

In this section, we collect several basic properties of the Tribonacci-Lucas numbers and present auxiliary results needed later. We also briefly recall tools from the theory of linear forms in three logarithms, as well as the Baker--Davenport reduction method. These ingredients play an essential role in the proof of our main theorem.

\subsection{Some properties of the Tribonacci-Lucas sequence}
We recall several essential properties of the Tribonacci-Lucas sequence. 
Recall that the characteristic equation 
\[
\varphi(x):=x^{3}-x^{2}-x-1=0
\]
has roots \(\alpha\), \(\beta\) and \(\gamma=\overline{\beta}\), in which
\begin{eqnarray}\label{Trib1}
	\alpha =\dfrac{1+(w_1+w_2)}{3}, \qquad 
	\beta = \dfrac{2-(w_1+w_2)+\sqrt{-3}\,(w_1-w_2)}{6},
\end{eqnarray}
with
\begin{eqnarray}\label{Trib2}
	w_1=\sqrt[3]{19+3\sqrt{33}}, \qquad 
	w_2=\sqrt[3]{19-3\sqrt{33}} .
\end{eqnarray}

The Binet's formula for the general term of the sequence can be written as
\begin{eqnarray}\label{Trib3}
	S_n = \alpha^{n}+\beta^{n}+\gamma^{n} \qquad (n\ge 0).
\end{eqnarray}
Numerically, the dominant root satisfies \begin{align}\label{Trib4}
	1.83 &< \alpha < 1.84    
\end{align} \noindent and the complex roots satisfy \begin{equation}\label{Tribk}
	0.73 < |\beta| = |\gamma| = \alpha^{-1/2} < 0.74 .
\end{equation}

From \eqref{Trib1}, \eqref{Trib2} and \eqref{Tribk} it is clear that the contribution of the complex conjugate roots \(\beta\) and \(\gamma\) to the right hand side of \eqref{Trib3} is very small.
In fact, if we set
\begin{equation}\label{Trib5}
	\xi(n) := S_n - \alpha^{n} = \beta^{n}+\gamma^{n},
\end{equation}
then
\[
|\xi(n)| \le \frac{2}{\alpha^{\,n/2}} \qquad (n\ge 1).
\]
The inequality follows directly from \(|\beta|=|\gamma|=\alpha^{-1/2}\); indeed for any \(n\ge 1\)
\[
|\xi(n)| = \bigl|\beta^{n}+\gamma^{n}\bigr| \le |\beta|^{n}+|\gamma|^{n}
= \alpha^{-n/2}+\alpha^{-n/2}
= 2\alpha^{-n/2} = \frac{2}{\alpha^{\,n/2}} .
\]
The following estimate taken from \cite{MD5} holds for every integer $m\geq 1.$

\begin{lemma}\label{lem1}
	For every integer \(m\ge 1\) we have
	\[
	\alpha^{m-1} \le S_m < \alpha^{m+1}.
	\]
\end{lemma}

Denote by \(\mathbb{K}:=\mathbb{Q}(\alpha,\beta)\) the splitting field of the polynomial \(\Psi\)
over \(\mathbb{Q}\). Then \([\mathbb{K}:\mathbb{Q}]=6\) and \([\mathbb{Q}(\alpha):\mathbb{Q}]=3\).
The Galois group of \(\mathbb{K}/\mathbb{Q}\) is
\[
\mathcal{G}:=\operatorname{Gal}(\mathbb{K}/\mathbb{Q})
\cong \{(1),\; (\alpha\beta),\; (\alpha\gamma),\;
(\beta\gamma),\; (\alpha\beta\gamma),\;
(\alpha\gamma\beta)\}.
\]
Thus we may identify each automorphism of \(\mathcal{G}\) with a permutation of the roots of \(\Psi\).
For instance, the permutation \((\alpha\beta)\) corresponds to the automorphism
\(\sigma : \alpha\mapsto\beta,\; \beta\mapsto\alpha,\; \gamma\mapsto\gamma\).
\subsection{Linear forms in logarithms}
The method we employ relies on the standard technique of establishing non zero lower bounds for specific linear forms in logarithms of algebraic numbers. Upper bounds, on the other hand, are derived by working with the associated Binet formula of the sequence. To obtain the required lower bounds we appeal to Baker's famous theorem on linear forms in logarithms. First we recall the definition of the logarithmic height of an algebraic number.

 Let \(\eta\) be an algebraic number of degree \(d\) whose minimal polynomial is
\[
P(x)=a_{0}\prod_{j=1}^{d}\bigl(x-\eta_{j}\bigr),
\]
where the leading coefficient \(a_{0}\) is positive and \(\eta_{1},\dots,\eta_{d}\) are the conjugates of \(\eta\).  The logarithmic height of \(\eta\) is defined as
\[
h(\eta):=\frac{1}{d}\Bigl(\log a_{0}+\sum_{j=1}^{d}
\log\bigl(\max\{|\eta_{j}|,1\}\bigr)\Bigr).
\]
If \(\eta=p/q\in\mathbb{Q}\) is a reduced rational number with \(q>0\), this definition simplifies to 
\(h(\eta)=\log\max\{|p|,q\}\).
The height function satisfies several useful properties, which we shall use freely in what follows:
\begin{equation}\label{heights}
	\begin{aligned}
		h(\eta_{1}\pm\eta_{2}) &\le h(\eta_{1})+h(\eta_{2})+\log 2,\\
		h(\eta_{1}\eta_{2}^{\pm1}) &\le h(\eta_{1})+h(\eta_{2}),\\
		h(\eta^{s}) &= |s|\,h(\eta) \qquad (s\in\mathbb{Z}).
	\end{aligned}
\end{equation}

For the lower bound we quote the version of Baker's theorem proved by Bugeaud, Mignotte and Siksek 
\cite{BMS} (Theorem 9.4), which itself goes back to Matveev.

\begin{theorem}[Bugeaud, Mignotte, Siksek, \cite{BMS}]\label{thm1}
	Let \(\eta_{1},\dots ,\eta_{t}\) be positive real algebraic numbers lying in a real algebraic number field
	\(\mathbb{K}\subset\mathbb{R}\) of degree \(D\).  
	Let \(b_{1},\dots ,b_{t}\) be nonzero integers such that
	\[
	\Gamma:=\eta_{1}^{b_{1}}\!\cdots\eta_{t}^{b_{t}}-1\neq 0 .
	\]
	Then
	\[
	\log|\Gamma| > -1.4\times 30^{\,t+3}\times t^{\,4.5}\times 
	D^{2}(1+\log D)(1+\log B)A_{1}\cdots A_{t},
	\]
	where
	\[
	B\ge \max\{|b_{1}|,\dots ,|b_{t}|\},
	\]
	and
	\[
	A_{j}\ge \max\bigl\{D\,h(\eta_{j}),\;|\!\log\eta_{j}|,\;0.16\bigr\},
	\qquad j=1,\dots ,t .
	\]
\end{theorem}

\subsection{Reduction procedure}

The bounds produced by Baker's theorem are typically far too large for direct computation.
To obtain sharper results we apply a Baker-Davenport type reduction.  The precise variant we need is due to Dujella and Peth\H{o} \cite{DP}.  
For a real number \(r\) we write \(\|r\|\) for the distance from \(r\) to the nearest integer, i.e.
\(\|r\|:=\min\{|r-n|:n\in\mathbb{Z}\}\).

\begin{lemma}[Dujella, Peth\H{o}, \cite{DP}]\label{DP}
	Let \(\kappa\neq0\), \(A\), \(B\) and \(\mu\) be real numbers with \(A>0\) and \(B>1\).  
	Let \(M>1\) be a positive integer and suppose that \(p/q\) is a convergent of the continued fraction
	expansion of \(\kappa\) satisfying \(q>6M\).  Put
	\[
	\varepsilon:=\|\mu q\|-M\|\kappa q\|.
	\]
	If \(\varepsilon>0\), then the inequality
	\[
	0<|m\kappa-n+\mu|<AB^{-k}
	\]
	has no solution in positive integers \(m,n,k\) with
	\[
	\frac{\log(Aq/\varepsilon)}{\log B}\le k \qquad\text{and}\qquad m\le M .
	\]
\end{lemma}

When \(\mu=0\) the quantity \(\varepsilon\) becomes negative, so Lemma \ref{DP} cannot be used directly.
In that case, we resort to a classical criterion of Legendre, a standard tool in Diophantine approximation.
The reader may consult the books of Cohen \cite{HC1,HC2} for further details.

\begin{lemma}[Legendre, \cite{HC1,HC2}]\label{lg}
	Let \(\kappa\) be a real number and \(x,y\) integers such that
	\[
	\Bigl|\kappa-\frac{x}{y}\Bigr|<\frac{1}{2y^{2}} .
	\]
	Then \(x/y\) is a convergent \(p_{k}/q_{k}\) of \(\kappa\).  
	Moreover, let \(M\) and \(N\) be non negative integers with \(q_{N}>M\) and set  
	\(a(M):=\max\{a_{i}:i=0,1,\dots ,N\}\), where the \(a_{i}\) are the partial quotients of \(\kappa\).
	The inequality
	\[
	\Bigl|\kappa-\frac{x}{y}\Bigr|\ge\frac{1}{(a(M)+2)\,y^{2}}
	\]
	holds for every pair \((x,y)\) of positive integers with \(0<y<M\).
\end{lemma}

Finally, we shall also make use of the following lemma due to G\'uzman S\'anchez and Luca
\cite{GSL} (Lemma 7).

\begin{lemma}[G\'{u}zman S\'{a}nchez, Luca, \cite{GSL}]\label{l3}
	Let \(r\ge1\) and \(H>0\) be numbers such that \(H>(4r^{2})^{r}\) and \(H>L/(\log L)^{r}\).  
	Then
	\[
	L<2^{r}H(\log H)^{r}.
	\]
\end{lemma}
\section{Proof of the Main Result}

\subsection{The case $ k\in [4,900] $}
One can see that the first four terms of the Tribonacci-Lucas sequence do not satisfy the Diophantine equation \eqref{main} since they are single digits. Hence for values of $ k\in [4,900] $, we write {\tt SageMath code} to check for all solutions to our Diophantine equation \eqref{main}. The only solutions obtained are the ones stated in Theorem \ref{thm1x}. Hence, from now onwards we assume $k>900.$

\subsection{The initial bound on $k$}
Assume \eqref{main} holds, let us examine the relationship between the variables in \eqref{main}. The number of digits of $S_n$ can be expressed as \begin{align*}
	d&=\lfloor\log_{10} {S_n} \rfloor +1 \leq \log_{10} {S_n} +1\leq \log_{10}{\alpha^{n+1}}+1 =(n+1)\log_{10} \alpha +1<n+2,
\end{align*} also
\begin{align*}
	d=\lfloor\log_{10} {S_n} \rfloor +1 >  \log_{10} {S_n} \geq \log_{10} {\alpha^{n-1}}=(n-1)\log \alpha >\frac{n-1}{2},
\end{align*} hence we have that \begin{align}\label{dbound}
	\frac{n-1}{2}<d<n+2.
\end{align}
\noindent Also from $$S_n<10^{\log_{10} {S_n}}<10^d<10^{\log _{10}{Sn}+1}=10S_n,$$ we have that \begin{align}
	S_n<10^d<10S_n.
\end{align}
\noindent From the above inequality and \eqref{lem1} we can see that 
\begin{align*}
	\alpha^{k-1}\leq S_k=S_m \cdot 10^d +S_n < 10S_mS_n +S_mS_n=11S_mS_n<\alpha^{n+m+2},
\end{align*} and  \begin{align*}
	\alpha^{k+1} \geq S_k =S_m\cdot 10^d +S_n>S_mS_n +S_n >S_mS_n> \alpha^{m+n-2},
\end{align*} thus we get \begin{align}
	m+n-3<k<m+n+3.
\end{align}
\noindent Also since $S_k=S_n \cdot 10^d +S_n > S_nS_m +S_n \geq 2S_n,$ we can see that $k\geq n+2$ and thus $k-n\geq 2.$    
Next, we prove the following lemma. 
\begin{lemma}\label{5}
	All solutions of the Diophantine equation \eqref{main} satisfy the inequality \begin{equation*}
		m<2.2\cdot 10^{11}\log k.
	\end{equation*}
\end{lemma}
\begin{proof}
	Suppose $k>900$, we rewrite \eqref{main} as \begin{align*}
		\alpha^k +\beta^k + \gamma^k =10^d (\alpha^m+\beta^m+\gamma^m)+S_n, \\ \alpha^k-10^d \alpha^m  =10^d (\beta^m+\gamma^m)+S_n-\beta^k-\gamma^k. 
	\end{align*}
	\noindent Dividing both sides by $10^d\alpha^m$ and taking absolute values on both sides yields \begin{align*}
		\left|\frac{\alpha^{k-m}}{10^d} -1\right|&\leq \left|\frac{\beta^m+\gamma^m}{\alpha^m} \right|+\left|\frac{S_n}{\alpha^m\cdot 10^d} \right|+\left|\frac{\beta^k}{\alpha^m \cdot 10^d} \right|+\left|\frac{\gamma^k}{\alpha^m\cdot 10^d} \right| \\ &< \frac{1}{\alpha^m} \left(\frac{2}{\alpha^{\frac{m}{2}}}+1+\frac{2}{\alpha^{\frac{k}{2}}\cdot 10^d } \right) < \frac{3}{\alpha^m},
	\end{align*} thus \begin{align}\label{lambda1}
		\left| \frac{\alpha^{k-m}}{10^d}-1 \right| < \frac{3}{\alpha^m},
	\end{align} where we have used the fact that $S_n<10^d, k>900,m\geq0, d\geq 1.$
	
	\noindent Fix
	\begin{equation}
		\Gamma_1 := \frac{\alpha^{k-m}}{10^d}-1 .
	\end{equation}
	
	\noindent Suppose $\Gamma_1=0,$ we get $\alpha^{k-m}=10^d$. This implies that $\alpha^{k-m}\in \mathbb{Q}$ which is impossible for all $k-m>0.$ Clearly we see that $\Gamma_1\neq 0.$
	Define 
	\begin{equation*}
		\eta_1 := \alpha, \ \eta_2 := 10,  \ b_1 := k-m, \ b_2 := -d, , \ t: = 2.
	\end{equation*}
	\noindent Since $d<n+2<(k-m)+5<k+5,$ we take $B=k+5.$
	\noindent Furthermore, 
	\begin{equation*}
		h(\eta_1) = h \left( \alpha \right) =\frac{1}{3} \log \alpha. 
	\end{equation*}
	
	\noindent We also have that $h(\eta_2) = h(10) = \log {10}. $ Hence, we let
	
	\begin{equation*}
		A_1 := \log \alpha, \ A_2 := 3\log {10}.
	\end{equation*}
	
	\noindent Thus, proceeding via Theorem \ref{thm1} yields 
	
	\begin{align*}
		\log |\Gamma_1| & > -1.4 \cdot 30^6 \cdot  2^{4.5} \cdot  3^2 (1 + \log 3)(1 + \log {(k+5)})(\log \alpha)(3\log {10}) \\&> -1.3\times 10^{11}\log k.
	\end{align*}
	where we have used the fact that 
	$$
	1+\log {(k+5)}=1+\log {k\left(1+\frac{5}{k}\right)}<1+\log {2k}=1+\log 2 +\log k <2\log k.
	$$
	\noindent Also \begin{equation*}
		\log|\Gamma_1| < \log 3 -m\log \alpha
	\end{equation*}
	Comparing the two bounds on $\log|\Gamma_1|$ yields 
	\begin{equation}
		m<2.2\cdot 10^{11}\log k. 
\end{equation} \end{proof}
We now prove the following lemma.
\begin{lemma}
	All solutions to the Diophantine equation \eqref{main} satisfy the inequality 
	\begin{align*}
		k<1.5\cdot 10^{29}.
	\end{align*}   
\end{lemma}
\begin{proof}
	We rewrite \eqref{main} as \begin{align*}
		\alpha^k +\beta^k+\gamma^k &=S_m \cdot 10^d + \alpha^n +\beta^n+\gamma^n, \\ \alpha^k -\alpha^n -S_m \cdot 10^d &= \beta^n +\gamma^n -(\beta^k+\gamma^n).
	\end{align*}    
	\noindent Factorizing by $\alpha^k$ yields; \begin{align*}
		\alpha^k \left(1-\alpha^{n-k}\right)-S_m\cdot 10^d=\beta^n +\gamma^n -(\beta^k+\gamma^n).
	\end{align*} 
	\noindent Dividing through by $\alpha^k \left(1-\alpha^{n-k}\right)$ and taking absolutes on both sides yields
	\begin{align*}
		\left|1-\frac{S_m\cdot 10^d}{\alpha^k (1-\alpha^{n-k})} \right| &\leq \left|\frac{\beta^n+\gamma^n}{\alpha^k \left(1-\alpha^{n-k}\right)} \right|+\left|\frac{\alpha^k +\beta^k}{\alpha^k \left(1-\alpha^{n-k}\right)} \right| \\ &< \frac{2}{\alpha^{\frac{n}{2}}\alpha^k \left(1-\alpha^{n-k}\right)}+\frac{2}{\alpha^{\frac{k}{2}}\alpha^k \left(1-\alpha^{n-k}\right)}, 
	\end{align*} \noindent hence 
	\begin{align*}
		\left|1-\frac{S_m\cdot 10^d}{\alpha^k (1-\alpha^{n-k})} \right| &< \frac{2}{\alpha^k} \left|\frac{1}{1-\alpha^{n-k}} \right|\left( \frac{1}{\alpha^{\frac{n}{2}}}+\frac{1}{\alpha^{\frac{k}{2}}}\right)<\frac{4}{\alpha^k},
	\end{align*} \noindent where we have used the fact that $n\geq 0, k>900$ and $k-n\geq2.$ Fix
	$$
	\Gamma_2: =1-\frac{S_m\cdot 10^d}{\alpha^k (1-\alpha^{n-k})}.
	$$
	Then, we have
	\begin{equation}\label{g2}
		|\Gamma_2| <\frac{4}{\alpha^k}.
	\end{equation} 
	Exactly, the same argument used in the proof of lemma \ref{5} shows that $\Gamma_2\neq 0$. Define  
	\begin{equation*}
		\eta_1: = \alpha , \ \eta_2: = 10, \ \eta_3: = \frac{S_m}{1-\alpha^{n-k}}, \ b_1: = -k, \ b_2 := d, \ b_3 := 1, \ t := 3.
	\end{equation*}
	We have $B:=k+5$, \begin{equation*}
		h(\eta_1):=h(\alpha)=\frac{1}{3}\log \alpha, \ h(\eta_2):=h(10)=\log {10},  
	\end{equation*} thus \begin{align*}
		A_1:= \log \alpha, \quad A_2:=3\log {10}.
	\end{align*}
	\begin{align*}
		h(\eta_3)&=h\left(\frac{S_m}{1-\alpha^{n-k}} \right) \\ &\leq h(S_m)+ \frac{k-n}{3} \log \alpha + \log 2 \\ &< \frac{m+1}{3} \log \alpha + \frac{m+3}{3} \log \alpha + \log 2 \\ &<9\cdot 10^{10}\log k ,  
	\end{align*} 
	where we have used the result in Lemma \ref{5}.
	Thus $A_3:= 2.7\cdot 10^{11}\log k.$ We proceed via Theorem \ref{thm1} and get 
	\begin{align*}
		\log |\Gamma_2| & > -1.4 \cdot  30^6 \cdot 3^{4.5} \cdot 3^2 \cdot (1 + \log 3)(1 + \log {(k+5)})(\log \alpha)(3 \log 10)(2.7\cdot 10^{11}\log k)\\
		& > -6.2 \cdot 10^{24}(\log k)^2.
	\end{align*}
	Also 
	\begin{equation}
		\log \Gamma_2 < \log 4 -k\log \alpha.   
	\end{equation}
	Comparing the two inequalities, we have that
	
	\begin{align}\label{eq7}
		k<1.1\cdot 10^{25}(\log k)^2.\end{align}
	Proceeding via Lemma \ref{l3} with parameters $H:=1.1\cdot 10^{25}, L:=k$ and $r=2$ yields
	\begin{equation}
		k<1.5\cdot 10^{29}.
	\end{equation} 
\end{proof}

\subsection{Reduction of the bounds}
The bounds produced by Theorem \ref{thm1} are still far too large for practical computation. We therefore need to reduce them substantially, which we accomplish by repeatedly applying Lemma~\ref{DP}.

First, we revisit \eqref{lambda1}and fix 
\begin{align*}
	\Lambda_1:= (k-m)\log \alpha -d\log {10}. \end{align*}
We already showed that $\Gamma_1\neq 0$ and thus $\Lambda_1\neq 0.$ Also, if $m\geq3$, we have $|\Gamma_1|=|e^{|\Lambda_1|}-1|<0.5$, which yields $e^{|\Lambda_1|}\leq 1+|\Gamma_1|<1.5.$ Thus \begin{align*}
	\left|d\log {10}-(k-m)\log \alpha \right|< \frac{4.5}{\alpha^m}   
\end{align*} and hence \begin{equation}
	0<\left|\frac{\log \alpha}{\log {10}}-\frac{d}{k-m} \right| <\frac{2}{(k-m)\alpha^m}.
\end{equation}
\noindent From Lemma \ref{lg}, we have that $\frac{d}{k-m}$ is a convergent to the continued fraction expansion of $\frac{\log \alpha}{\log {10}}$, say $\frac{p_i}{q_i}$. Since $q_i\leq k-m<k<1.5\cdot 10^{29}$, we obtain the upper bound for indices $i=i_0$ and thus the maximum value of $a_{\text{max}}=\text{max}[a_0, a_1, \cdots a_{i_0}]$ of the continued fraction of $\frac{\log \alpha}{\log {10}}$. By applying a property of continued fraction, we write \begin{align}\label{m1}
	\frac{1}{(a_{\text{max}}+2)(k-m)^2} \leq \frac{1}{(a_i+2)(k-m)^2}<\left|\frac{\log \alpha}{\log {10}}-\frac{d}{k-m}\right|<\frac{2}{(k-m)\alpha^m} .
\end{align}
A simple computer program in Sage-math yields 
$$
[a_1,a_1,a_2,\cdots ]=[0, 3, 1, 3, 1, 1, 14, 1, 3, 3, 6, 1, 13, 3, \cdots ],
$$ 
for which we obtain $a_{\text{max}}=a_{34}=44$. Thus from \eqref{m1} we get \begin{align}
	m<\frac{\log{\left( 2\times46\times 1.5\cdot 10^{29}\right)}}{\log {\alpha}}<119.
\end{align}
With the new bound on $m$, we obtain a new bound on $k$ which is \begin{align}
	k<6.3\cdot10^{14}.
\end{align}
\noindent Now let \begin{align}
	\Lambda_2:= d\log {10}-k\log \alpha +\log {\left(\frac{S_m}{1-\alpha^{n-k}} \right)}.
\end{align} \noindent From inequality \eqref{g2}, we have that \begin{align*}
	|\Gamma_2|=\left|e^{|\Lambda_2|} -1\right|<\frac{4}{\alpha^k}.
\end{align*} 
\noindent Note that if $k\geq 2$, then the right hand side of the above inequality is about \(\frac{4}{\alpha^2}<\frac{1}{2}\). Thus we get \begin{align*}
	|\Lambda_2|<\frac{6}{\alpha^k},
\end{align*} \noindent which implies that \begin{align*}
	\left|d\log {10}-k\log \alpha+\log {\left(\frac{S_m}{1-\alpha^{n-k}} \right)} \right|<\frac{6}{\alpha^k}.
\end{align*}
\noindent Dividing through by $\log \alpha$ gives \begin{align*}
	\left|d\frac{\log {10}}{\log \alpha}-k+\frac{\log {\left(\frac{S_m}{1-\alpha^{n-k}} \right)}}{\log \alpha} \right|<\frac{6}{\alpha^k\log \alpha},
\end{align*} \noindent thus we apply Lemma \ref{DP} with the quantities 
\begin{align*}
	\kappa:=\dfrac{\log 10}{\log \alpha}, \quad \mu:=\frac{\log {\left(\frac{S_m}{1-\alpha^{n-k}} \right)}}{\log \alpha}, \quad A:=\frac{6}{\log \alpha}, \quad B:=\alpha, \quad M:=6.3\cdot 10^{14}.
\end{align*}
\noindent With the help of a simple code in SageMath for $m\in [0,119], 1\leq k-n\leq m+3$ and $d<k+5<M$, we find that the convergent 
\begin{align*}
	\dfrac{p}{p}=\dfrac{p_{46}}{q_{46}}=\dfrac{3545917794526275601837767}{938425170962281070635339},
\end{align*}
is such that $ q=q_{46}>6M $. Furthermore, it gives $ \varepsilon > 0.0000000077243415 $, and thus,
\begin{align*}
	k\le \dfrac{\log((6/\log\alpha)q/\varepsilon)}{\log \alpha}<124.9872372.
\end{align*}
which yields to $k\leq 125.$ This contradicts our assumption that $k> 900$. Hence, Theorem \ref{thm1x} is completely proved. \qed

\section{Conclusion}
We determined all Tribonacci-Lucas numbers that can be expressed as concatenations of two terms of the same sequence, showing that only two such representations exist. The proof combines Baker’s theory on linear forms in logarithms with reduction techniques and computational verification, and these methods may be extended to similar problems for other recurrence sequences.

\section*{Acknowledgements} 

I thank the anonymous reviewers for their helpful comments and suggestions that improved the clarity of this work.

\makeatletter
\renewcommand{\@biblabel}[1]{[#1]\hfill}
\makeatother

\end{document}